\documentclass[11pt, reqno]{amsart}
\usepackage{amsmath, amsthm, amscd, amsfonts, amssymb, graphicx, color}
\usepackage[bookmarksnumbered, colorlinks, plainpages]{hyperref}
\numberwithin{equation}{section}
\DeclareMathOperator*{\essinf}{ess\,inf}
\DeclareMathOperator*{\esssup}{ess\,sup}

\newtheorem{theo}{Theorem}[section]
\newtheorem{defi}[theo]{Definition}

\newtheorem{prop}[theo]{Proposition}

\newtheorem{coro}[theo]{Corollary}
\newtheorem{lem}[theo]{Lemma}
\newtheorem{rem}[theo]{Remark}

\renewenvironment{proof}[1][Proof]{\textbf{#1. } }{\ \rule{0.5em}{0.5em}}

\begin{document}
\setcounter{page}{1}

\title[Calder\'{o}n's complex interpolations of variable-exponent Fofana spaces] {Calder\'{o}n's complex interpolations of variable-exponent Fofana spaces} 

\author[N. Diarra et al]{ Nouffou Diarra$^{1*}$, Pokou Nagacy$^2$ and Sekou Coulibaly$^3$} 

\maketitle
\address{$^{1}$ Laboratoire de Math\'{e}matiques et Applications, UFR Math\'{e}matiques et Informatique, Universit\'{e} F\'{e}lix Houphou\"{e}t-Boigny, 22 BP 582 Abidjan 22, C\^ote d'Ivoire.}\\
e-mail : nouffoud$@$yahoo.fr  

\address{ $^{2}$ Laboratoire des Sciences et Technologies de l'Environnement, UFR Environnement, Universit\'{e}
	Jean Lorougnon Gu\'{e}d\'{e}, BP 150 Daloa, C\^{o}te d'Ivoire.}\\
e-mail : pokounagacy@yahoo.com 

\address{ $^{3}$ D.E.R Math\'{e}matiques et Informatique, Universit\'{e} de Bamako, BPE 3206 Bamako, Mali}\\
e-mail : csekou91@yahoo.fr\\

\address{\textbf{*Corresponding author}: \textnormal{nouffoud$@$yahoo.fr}}\\
%
%








\begin{abstract}
Fofana spaces, introduced in the framework of Wiener amalgam spaces, can be viewed as generalizations of Morrey spaces. These spaces have attracted considerable attention in recent years and have found applications in the study of partial differential equations. In this paper, we investigate their variable-exponent counterparts. More specifically, we describe Calder\'{o}n's, first and second, complex interpolations of  Fofana spaces with variable exponent, as well as the corresponding interpolation spaces of their preduals. Our results recover, as special cases, those already known for (classical) Fofana and Morrey spaces. 
\end{abstract} 


\section{Introduction} \label{Sec1}
 
The interpolation theory has attracted substantial interest in modern functional analysis and has widely studied in many Banach spaces, including Morrey spaces. We recall that Morrey spaces are generalizations of Lebesgue spaces and have been introduced since 1938 by C. Morrey  \cite{Morrey} in order to study the regularity problem of solutions  for partial differential equations. In recent decades,
 a lot of  papers have been devoted to complex interpolations of  Morrey-type spaces (see, for example, \cite{Dia-Fof-2024,Hakim,Hakim-Sawano-2016,Hakim-saw-2016,Yang-al-2013}). 

In this paper, we are interested in the description of Calder\'{o}n's  complex interpolations of Fofana spaces with variable exponent. Note that classical Fofana spaces, introduced by I. Fofana \cite{Fof-1988} in 1988, consist of a class of function spaces which provide a natural setting for dilation and translation properties, combine local and global integrability conditions, and therefore exhibit a richer structure than Lebesgue and Morrey spaces. We refer to \cite{Coul-Fof-2019,Dia-Fof-2021,Dia-Fof-2023,Dosso-Fof-San,Feuto-BP,Fof-1988,Wei} and the references therein for a historical background of classical Fofana spaces. Let us mention that, Fofana spaces  can be viewed as genaralizations of Morrey spaces (see \cite{Dia-Fof-2021} for more precision).
Thus, a precise description of complex interpolation of Fofana spaces is not merely an abstract interpolation problem; but it offers a  framework for identifying intermediate spaces and extending classical boundedness results of operators to a broader families of Morrey-type function spaces. 

 Let $d$ be a fixed positive integer. For any measurable subset $\Omega$ of $\mathbb{R}^d$, $|\Omega|$ and $\chi_{\Omega}$ stand for its Lebesgue measure and its characteristic function, respectively.
Let  $p(\cdot)$ be a variable exponent on $\mathbb{R}^d$, that is, $p$ is a measurable function  from $\mathbb{R}^d$ to $[1,\infty]$. 
The conjuguate exponent $p^{\prime}(\cdot)$ of  $p(\cdot)$ is defined by the formula $\frac{1}{p^{\prime}(\cdot)}=1-\frac{1}{p(\cdot)}$, with the convention $\frac{1}{\infty}=0$.
We denote by $\mathcal{P}(\mathbb{R}^d)$ the set of all variable exponents $p(\cdot)$ on $\mathbb{R}^d$ such that $1<p_{-} \leq p(\cdot)\leq p_{+}<\infty$, where
\begin{align*}
	p_{-}=\essinf_{x \in  \mathbb{R}^d} p(x)  \mbox{ and } p_{+}=\esssup_{x \in  \mathbb{R}^d} p(x).
\end{align*}
The variable exponent Lebesgue space $L^{p(\cdot)}:=L^{p(\cdot)}(\mathbb{R}^d)$ consists of all measurable functions $f$ such that, for some $\lambda>0$,
\begin{align*}
	\int_{\mathbb{R}^d}\left( \frac{|f(x)|}{\lambda }\right)^{p(x)}dx<\infty.
\end{align*}
The space $L^{p(\cdot)}$ is a Banach space, when equipped with the Luxemburg norm
\begin{align*}
	\| f\|_{p(\cdot)}:=\inf \left\{ \lambda>0 :\int_{\mathbb{R}^d}\left( \frac{|f(x)|}{\lambda }\right)^{p(x)}dx \leq 1 \right\}.
\end{align*}
If  $p(\cdot)$ is a constant function  $p(x)=p$, $x \in \mathbb{R}^{d}$, then $\| \cdot\|_{p(\cdot)}$ coincides with the classical  Lebesgue norm $\| \cdot\|_{p}$ and so $L^{p(\cdot)}$ is equal to the classical Lebesgue space $L^{p}$.  
The space  $L^{p(\cdot)}_{\rm loc}$ is defined by
\begin{align*}
	L^{p(\cdot)}_{\rm loc}=\left\{ f : f\chi_{K}\in L^{p(\cdot)}, \mbox{ for all compact subset $K$ of $\mathbb{R}^d$}\right\}.
\end{align*} 
Let $f \in L^{p(\cdot)}_{\rm loc}$, $1\leq q \leq \infty$ and $r>0$.  We set 
\begin{align}
	_{r}\|f\|_{p(\cdot),q}:=\left \| \left\{ \|f\chi_{I_{k}^{r}}\|_{p(\cdot)} \right\}_{k\in \mathbb{Z}^d}\right\|_{\ell^{q}}=
	\left\{
	\begin{array}{rl}
		&\displaystyle \left(\sum_{k\in \mathbb{Z}^d} \|f\chi_{I_{k}^{r}}\|^{q}_{p(\cdot)} \right)^{\frac{1}{q}} \qquad \mbox{if } q<\infty ,\\
		&\displaystyle\sup_{k\in \mathbb{Z}^d} \|f\chi_{I_{k}^{r}}\|_{p(\cdot)} \qquad \qquad \, \mbox{ if }\ q =\infty,
	\end{array}
	\right.
\end{align}
where
\begin{align*}
	I_{k}^{r}=\displaystyle \prod_{j=1}^{d}\left[ rk_{j},r(k_{j}+1) \right), \quad k=(k_{1},k_{2},...,k_{d})\in \mathbb{Z}^{d}.
\end{align*}

\begin{defi}
	Let $p(\cdot)$ be a variable exponent on $\mathbb{R}^d$ and $ 1 \leq \alpha , q \leq \infty$. 
	\begin{enumerate}
	\item The variable exponent  amalgam space $\left(L^{p(\cdot)},\ell^{q}\right):=(L^{p(\cdot)}, \ell^{q})(\mathbb{R}^d)$ is defined by  
	\begin{align*}
		(L^{p(\cdot)}, \ell^{q})= \left\{ f\in L^{p(\cdot)}_{\rm{loc}} \quad : \quad 	_{1}\|f\|_{p(\cdot),q}<\infty \right\}.
	\end{align*}
	\item The variable-exponent  Fofana space $\left(L^{p(\cdot)},\ell^{q}\right)^{\alpha}:=\left(L^{p(\cdot)}, \ell^{q}\right)^{\alpha}(\mathbb{R}^d)$ is defined by 
	\begin{align*}
		\left(L^{p(\cdot)}, \ell^{q}\right)^{\alpha}= \left\{ f\in L^{p(\cdot)}_{\rm{loc}}\;:\;\|f\|_{p(\cdot),q,\alpha}<\infty \right\},
	\end{align*} 
	where
	\begin{align}\label{eq1.2}
		||f||_{p(\cdot),q,\alpha}= \sup_{r >0}\; r^{d\left(\frac{1}{\alpha}-\frac{1}{N_{p,r}}\right)}\; {_{r}||f||_{p(\cdot),q}}
	\end{align}
	with \begin{align*}
		N_{p,r} = \left\{
		\begin{array}{rl}
			&p_{+} \quad \mbox{ if } r \leq 1,\\
			&p_{-} \quad \mbox{ if } r > 1.
		\end{array}
		\right.
	\end{align*}
\end{enumerate}
\end{defi}
Variable-exponent Fofana spaces were introduced by Yang and Zhou \cite{yang-zhou-23} in 2023, where they proved that $\left(L^{p(\cdot)},\ell^{q}\right)^{\alpha}$ is non trivial if and only if $ p(\cdot) \leq \alpha \leq q$. Thus, we shall always assume that this condition is fulfilled. Note that, the norm in \eqref{eq1.2} can be replaced by
 \begin{align}
 	||f||_{p(\cdot),q,\alpha}= \sup_{r >0}\; \left \|St_{r}^{(\alpha)}f \right\|_{p(\cdot),q},
 \end{align}
 where the dilation operator $St_{r}^{(\alpha)}$ is defined by
 \begin{align*}
 	St_{r}^{(\alpha)}f(x)= r^{-\frac{d}{\alpha}}f\left(r^{-1} x \right), \quad x\in \mathbb{R}^{d}.
 \end{align*}
 Moreover, if $p(\cdot)=p$ is a constant exponent then $\left(L^{p(\cdot)},\ell^{q}\right)^{\alpha}$ coincides with the  classical Fofana space $\left(L^{p},\ell^{q}\right)^{\alpha}$. 
 
 Now, let us recall the definition of predual spaces of variable-exponent Fofana spaces.
 
\begin{defi}
	Let $p(\cdot)$ be a variable exponent on $\mathbb{R}^d$ and $ 1 \leq q \leq \alpha \leq p(\cdot) \leq \infty$.
	\begin{enumerate} 
		\item  A sequence $\left\{(c_{n},\rho_{n},f_{n})\right\}_{n\geq 1}$ of elements of $\;\mathbb{C}\times (0, \infty)\times (L^{p(\cdot)}, \ell^{q})$ is called a $\mathfrak{h}$-decomposition of a function $f$  if $f=\displaystyle\sum_{n\geq 1} c_{n}St_{\rho_{n}}^{(\alpha)}f_{n}$ in the sense of $L^{p(\cdot)}_{\rm{loc}}$ and
		\begin{align*}
			\left\{
			\begin{array}{rl}
				&\|f_{n}\|_{p(\cdot),q}\leq 1, \quad n\geq 1,\\\\
				&\displaystyle \sum_{n\geq 1}|c_{n}|< \infty.
			\end{array}
			\right.
		\end{align*}
		\item The space $\mathcal{H}(p(\cdot),q,\alpha):=\mathcal{H}(p(\cdot),q,\alpha)(\mathbb{R}^d)$ is defined as the set of all $f \in L^{p(\cdot)}_{\rm{loc}}$  whose set of  $\mathfrak{h}$-decompositions is nonvoid.
		\item For any $f\in \mathcal{H}(p(\cdot),q,\alpha)$, we set 
		\begin{align*}
			\|f\|_{\mathcal{H}(p(\cdot),q,\alpha)}=\inf \left\{\displaystyle \sum_{n\geq 1}|c_{n}| : f=\displaystyle\sum_{n\geq 1} c_{n}St_{\rho_{n}}^{(\alpha')}f_{n}\right\}.
		\end{align*}
\end{enumerate}
\end{defi}

	Let $X$ be a Banach function space on $\mathbb{R}^d$, which norm is denoted by $\|\cdot\|_{X}$. We use here the terminology from the theory of Banach function spaces in the sense of Zaanen  (see \cite[Chapter 15]{Zaanen}).
	The space $X$ is said to satisfy the Fatou property if, for any sequence $(f_{n})_{n\geq 1}$ of elements of $X$, we have 
		\begin{align*}
			\left[0\leq f_{n} \uparrow f \quad \mbox{ and } \displaystyle \sup_{n\geq 1}\|f_{n}\|_{X}<\infty \right]\Longrightarrow \left[ f\in X\;\mbox{ and } \displaystyle \lim_{n\rightarrow \infty} \|f_{n}\|_{X}=\|f\|_{X}\right].
		\end{align*}
	 The space $X$ is said to have an absolutely continuous norm if, for any $f \in X$, $\displaystyle \lim_{n\rightarrow \infty}\|f\chi_{E_{n}}\|_{X}=0$ whenever  $(E_{n})_{n\geq 1}$ is a nonincreasing 
	 sequence   of  measurable subsets of $\mathbb{R}^{d}$ satisfying $\left|\displaystyle \bigcap_{n\geq 1}E_{n}\right|=0$. We denote by $X^{*}$  and $X^{\prime}$ the topological dual and the K\"{o}the dual (or associate space), respectively,   of $X$.
	  
	We summarize below some useful properties of variable-exponent Fofana spaces and their preduals (see \cite{Dia-2023,Dia-Nag-Kpa}).
	 
	 \begin{prop} \label{P2.9}
	 	Let $p(\cdot)$ be a variable exponent on $\mathbb{R}^d$ and $ p(\cdot) \leq  \alpha \leq q \leq \infty$. Then, 
	 	\begin{enumerate}
	 		\item the space $\left(L^{p(\cdot)},\ell^{q}\right)^{\alpha}$ is a Banach function space satisfying the Fatou
	 		property;
	 		\item if  $p(\cdot) \in \mathcal{P}(\mathbb{R}^d)$ and $\alpha<\infty$
	 		then the subset $\left(L^{p(\cdot)},\ell^{q}\right)^{\alpha}_{\rm a}$ of $\left(L^{p(\cdot)},\ell^{q}\right)^{\alpha}$ having an absolutely continuous norm is the closure of the set $\mathcal{C}_{\rm{c}}^{\infty}$ of all infinitely differentiable and compactly supported functions on $\mathbb{R}^{d}$ in $\left(L^{p(\cdot)},\ell^{q}\right)^{\alpha}$;
	 		\item $\mathcal{H}(p^{\prime}(\cdot),q^{\prime},\alpha^{\prime})$ is a Banach function space;  
	 		\item if  $p(\cdot) \in \mathcal{P}(\mathbb{R}^d)$ then $\mathcal{H}(p^{\prime}(\cdot),q^{\prime},\alpha^{\prime})$ has an absolutely continuous norm,  the set $\mathcal{C}_{\rm{c}}^{\infty}$ is dense in $\mathcal{H}(p^{\prime}(\cdot),q^{\prime},\alpha^{\prime})$ and
	 		\begin{align}\label{eq2.4}
	 			\left[\mathcal{H}(p^{\prime}(\cdot),q^{\prime},\alpha^{\prime})\right]^{*}=\left[\mathcal{H}(p^{\prime}(\cdot),q^{\prime},\alpha^{\prime})\right]^{\prime}=\left(L^{p(\cdot)},\ell^{q}\right)^{\alpha};
	 		\end{align}
	 		\item  if  $p(\cdot) \in \mathcal{P}(\mathbb{R}^d)$ and $ p(\cdot) < \alpha < q \leq \infty$ then $\mathcal{H}(p^{\prime}(\cdot),q^{\prime},\alpha^{\prime})$ satisfies the Fatou property,
	 		\begin{align}\label{eq2.5}
	 			\left[\left(L^{p(\cdot)},\ell^{q}\right)^{\alpha}\right]^{\prime}=\mathcal{H}(p^{\prime}(\cdot),q^{\prime},\alpha^{\prime})
	 		\end{align}
	 		and 
	 		\begin{align}\label{eq2.6}
	 			\left[\left(L^{p(\cdot)},\ell^{q}\right)^{\alpha}_{\rm a}\right]^{*}=\mathcal{H}(p^{\prime}(\cdot),q^{\prime},\alpha^{\prime}).
	 		\end{align}
	 	\end{enumerate}
	 \end{prop}
	 
	 The main purpose of this paper is to describe Calder\'{o}n's, first and second, complex interpolation spaces of the  variable-exponent Fofana spaces $\left(L^{p(\cdot)},\ell^{q}\right)^{\alpha}$ and the corresponding interpolation spaces of their preduals $\mathcal{H}(p^{\prime}(\cdot),q^{\prime},\alpha^{\prime})$. Note that there are some recent results on the description of complex interpolations of  classical Fofana spaces (see \cite{Dia-Fof-2024}), which are special cases of the main results established in the  present paper.
	 
The remainder of this paper is structured as follows. In Section \ref{Sec2}, we recall some useful properties of Calder\'{o}n's complex interpolation methods. Section \ref{Sec3} deals with some auxiliary results on variable-exponent Fofana spaces that will play an important role in establishing the main results. 
Section \ref{Sec4}  is devoted to presenting complex interpolation results on variable-exponent Fofana spaces, as well as their preduals.


\section{Background on Calder\'{o}n's complex interpolation methods} \label{Sec2}

 A couple of Banach spaces $\left(X_{0},X_{1}\right)$ is said to be compatible if there exists a Hausdorff topological vector space in which $X_{0}$ and $X_{1}$ are continuously embedded.
Let us recall the definition of Calder\'{o}n's product. 
\begin{defi}
	Let $\theta \in (0,1)$ and $\left(X_{0},X_{1}\right)$  be a compatible couple of Banach function spaces.
	The Calder\'{o}n product $X_{0}^{1-\theta}X_{1}^{\theta}$ is defined by
	\begin{align*}
		X_{0}^{1-\theta}X_{1}^{\theta}=\bigcup_{f_{0}\in X_{0},\,f_{1}\in X_{1}}\left\{ f:\mathbb{R}^{d}\rightarrow \mathbb{C}\;:\;|f|\leq |f_{0}|^{1-\theta}|f_{1}|^{\theta}\;\;\mbox{in $L^{0}$}\right\}
	\end{align*}
	and equipped with the norm
	\begin{align*}
		\|f\|_{X_{0}^{1-\theta}X_{1}^{\theta}}:=\inf \left\{\|f_{0}\|_{X_{0}}^{1-\theta}\|f_{1}\|_{X_{1}}^{\theta} : f_{0}\in X_{0},f_{1}\in X_{1} \mbox{ and }|f|\leq|f_{0}|^{1-\theta}|f_{1}|^{\theta} \mbox{ in } L^{0} \right\}.
	\end{align*}
\end{defi}

In this paper, we focus on Calder\'{o}n's complex interpolation functors $\left( X_{0},X_{1}\right)\mapsto [X_{0},X_{1}]_{\theta}$  and $\left( X_{0},X_{1}\right)\mapsto [X_{0},X_{1}]^{\theta}$, introduced by A. P. Calder\'{o}n  \cite{Calderon}. For their definitions, we refer the reader to \cite{Berg-löfs,Calderon}.
For the proofs of our main results, we shall need the following well known properties.

\begin{prop}\label{prop 2.5}
	Let $\theta \in (0,1)$ and $\left(X_{0},X_{1}\right)$ be a compatible couple of Banach function spaces.
	\begin{enumerate}
	\item \cite{W.Yuan-2016}. We have
	\begin{align}\label{eq 2.1}
		\left[ X_{0},X_{1}\right]_{\theta}=\overline{X_{0}\cap X_{1}}^{X_{0}^{1-\theta}X_{1}^{\theta}}.
	\end{align}
	\item  \cite[Theorem 1.14, chapter IV]{Kr-Pe-Se}. If $X_{0}^{1-\theta}X_{1}^{\theta}\;$ has an absolutely continuous norm, then
	\begin{align}\label{eq 2.2}
		\left[ X_{0},X_{1}\right]_{\theta}=X_{0}^{1-\theta}X_{1}^{\theta}.
	\end{align}
	\item \cite[Theorem 3.4]{Mastylo-saw}. If the unit ball of $X_{0}^{1-\theta}X_{1}^{\theta}$ is closed in $X_{0}+X_{1}$, then 
	\begin{align}\label{eq 2.3}
		\left[X_{0},X_{1}\right]^{\theta}=X_{0}^{1-\theta}X_{1}^{\theta}.
	\end{align}
\end{enumerate}
\end{prop}

\begin{rem}\label{rem 2.6}
	It follows from the remark just after \cite[Theorem 1.14, chapter IV]{Kr-Pe-Se} that if $X_{0}$ or $X_{1}$ has an absolutely continuous norm, then $X_{0}^{1-\theta}X_{1}^{\theta}$ has also an absolutely continuous norm and so \eqref{eq 2.2} holds. 
\end{rem}

We have the following result concerning the relationship between the Calder\'{o}n product of Banach function spaces and their K\"{o}the dual spaces.

\begin{prop}  (\cite[Theorem 2.10 ]{Schep})\label{prop 2.7}.
	Let $\theta \in (0,1)$. If $X_{0}$ and $X_{1}$ are two Banach function spaces  on $\mathbb{R}^{d}$ satisfying the Fatou property, then 
	\begin{align*}
		\left(X_{0}^{1-\theta}X_{1}^{\theta}\right)^{\prime}=(X_{0}^{\prime})^{1-\theta}(X_{1}^{\prime})^{\theta}.
	\end{align*}
\end{prop} 

The following result is a particular case of \cite[Theorem 4.5.1]{Berg-löfs}.
\begin{prop}\label{prop 2.8}
Let $\theta \in (0,1)$ and $\left(X_{0},X_{1}\right)$ be a compatible couple of Banach spaces.
	If  $X_{0}\cap X_{1}$  is dense in both $X_{0}$ and  $X_{1}$, then
	\begin{align*}
		\left(\left[X_{0},X_{1}\right]_{\theta}\right)^{*}=\left[X_{0}^{*},X_{1}^{*}\right]^{\theta}.
	\end{align*}
\end{prop}




 \section{Auxiliary results on variable-exponent Fofana spaces} \label{Sec3}
 In this subsection, we establish some properties that will play a key role in establishing our main results. Let us begin with the following definitions. 
 \begin{defi}
 		Let $p(\cdot)$ be a variable exponent on $\mathbb{R}^d$ and $ p(\cdot) \leq  \alpha \leq q \leq \infty$.
 \begin{enumerate}
 	\item The space $V^{\infty}(L^{p(\cdot)}, \ell^{q})^{\alpha}$ is defined as the closure of $L^{\infty}\cap(L^{p(\cdot)}, \ell^{q})^{\alpha}$ in $(L^{p(\cdot)}, \ell^{q})^{\alpha}$.
 	\item  The space $V^{0}(L^{p(\cdot)}, \ell^{q})^{\alpha}$ is defined as follows:
 		\begin{align*}
 		V^{0}(L^{p(\cdot)}, \ell^{q})^{\alpha}= \left\{ f\in (L^{p(\cdot)}, \ell^{q})^{\alpha}: \displaystyle \lim_{m\rightarrow 0} \left\| \chi_{[0,m]}(|f|)f \right\|_{p(\cdot),q,\alpha}=0 \right\}.
 	\end{align*}
 \end{enumerate}
 \end{defi}

We have the following result.
 \begin{prop}
 	Let $p(\cdot)$ be a variable exponent on $\mathbb{R}^d$ and $ p(\cdot) \leq  \alpha \leq q \leq \infty$. Then 
 	\begin{align*}
 		V^{\infty}(L^{p(\cdot)}, \ell^{q})^{\alpha}= \left\{ f\in (L^{p(\cdot)}, \ell^{q})^{\alpha}: \displaystyle \lim_{m\rightarrow \infty} \left\| \chi_{\left(m,\infty\right)}(|f|)f \right\|_{p(\cdot),q,\alpha}=0 \right\}.
 	\end{align*}
 \end{prop}
 
 \begin{proof}
 	Let $f$ be an element of $(L^{p(\cdot)}, \ell^{q})^{\alpha}$.\\
 	$\bullet$ Let $\epsilon$ and $m$ be two positive real numbers and suppose that $f$ is in $V^{\infty}(L^{p(\cdot)}, \ell^{q})^{\alpha}$. Then, there exists $g\in L^{\infty}\cap(L^{p(\cdot)}, \ell^{q})^{\alpha}$ such that $\|f-g\|_{p(\cdot),q,\alpha}<\epsilon$. Note that 
 	\begin{align*}
 		\left|f\chi_{\{|f| > m\}}\right|
 		&\leq  \left|f\chi_{\{|f| > m\}\cap\{|g| \leq \frac{m}{2}\}}\right| +  \left|f\chi_{\{|f| > m\}\cap\{|g| > \frac{m}{2}\} }\right| \\
 		&\leq  \left|f\chi_{\{|f| > m\}\cap\{|g| \leq \frac{m}{2}\}}\right|+  \left|(f-g) \chi_{\{|f| > m\}\cap\{|g| > \frac{m}{2}\} }\right| +  \left|g\chi_{\{|f| > m\}\cap\{|g| > \frac{m}{2}\} }\right| \\
 		&\leq 3 |f-g|+ \left|g\chi_{\{|g| > \frac{m}{2}\} }\right|
 	\end{align*}
 	and therefore, for $m\geq 2\|g\|_{\infty}$, we get 
 	\begin{align*}
 		\|f\chi_{\{|f| > m\}}\|_{p(\cdot),q,\alpha} \leq  3\|f-g\|_{p(\cdot),q,\alpha} <3\epsilon.
 	\end{align*}
 	Consequently, we obtain $\displaystyle \lim_{m\rightarrow \infty} \left\| \chi_{\left(m,\infty\right)}(|f|)f \right\|_{p(\cdot),q,\alpha}=0 $. Thus 
 	\begin{align} \label{eq3.1}
 	V^{\infty}(L^{p(\cdot)}, \ell^{q})^{\alpha} \subset \left\{ f\in (L^{p(\cdot)}, \ell^{q})^{\alpha}: \displaystyle \lim_{m\rightarrow \infty} \left\| \chi_{\left(m,\infty\right)}(|f|)f \right\|_{p(\cdot),q,\alpha}=0 \right\}.
 \end{align}
 	$\bullet$ Conversely, suppose that $\displaystyle \lim_{m\rightarrow \infty} \left\| \chi_{\left(m,\infty\right)}(|f|)f \right\|_{p(\cdot),q,\alpha}=0 $. Let us set, for any integer $n\geq 1$, $f_{n}=f\chi_{\{|f|\leq n\}}$. We have
 	\[
 	\left\{
 	\begin{array}{rl}
 		&f_{n} \in L^{\infty}\cap(L^{p(\cdot)}, \ell^{q})^{\alpha} \;,\;\;\;\;\;\;\;\;\;\;n\geq 1 ,\\\\
 		&\displaystyle \lim_{n\rightarrow \infty} \left\| f_{n}-f \right\|_{p(\cdot),q,\alpha} = \displaystyle \lim_{n\rightarrow \infty} \left\| f\chi_{\{|f|> n\}} \right\|_{p(\cdot),q,\alpha}=0.
 	\end{array}
 	\right.
 	\]
 	Consequently, $f \in V^{\infty}(L^{p(\cdot)}, \ell^{q})^{\alpha}$. Thus
 	\begin{align}\label{eq3.2}
 		\left\{ f\in (L^{p(\cdot)}, \ell^{q})^{\alpha}: \displaystyle \lim_{m\rightarrow \infty} \left\| \chi_{\left(m,\infty\right)}(|f|)f \right\|_{p(\cdot),q,\alpha}=0 \right\} \subset V^{\infty}(L^{p(\cdot)}, \ell^{q})^{\alpha}.
 	\end{align}
 By	combining \eqref{eq3.1} and \eqref{eq3.2}, we obtain the desired result.
 \end{proof}
 
The following result holds.
 \begin{prop}
 		Let $p(\cdot)$ be a variable exponent on $\mathbb{R}^d$ and $ p(\cdot) \leq  \alpha \leq q \leq \infty$. Then $V^{0}(L^{p(\cdot)}, \ell^{q})^{\alpha}$ is a closed subset of  $(L^{p(\cdot)}, \ell^{q})^{\alpha}$.
 \end{prop}
 
 \begin{proof}
 	Suppose that $(f_{n})_{n\geq 1}$ is a sequence of elements of $V^{0}(L^{p(\cdot)}, \ell^{q})^{\alpha}$ which converges in $(L^{p(\cdot)}, \ell^{q})^{\alpha}$ to an element $f$. Let us consider a positive real number $m$ and an integer $n\geq 1$. We have 
 	\begin{align*}
 	|f|\chi_{\{|f| \leq m\}}
 	&\leq |f-f_{n}|\chi_{\{|f| \leq m\}}+ |f_{n}|\chi_{\{|f| \leq m\}\cap\{|f_{n}| \leq 2m\} }+ |f_{n}|\chi_{\{|f| \leq m\}\cap\{|f_{n}| > 2m\}} \\
 	&\leq |f-f_{n}|+ |f_{n}|\chi_{\{|f_{n}| \leq 2m\}} +2|f_{n}-f|= 3 |f-f_{n}|+ |f_{n}|\chi_{\{|f_{n}| \leq 2m\}}
 \end{align*}
 and therefore
 \begin{align*}
 	\left\| f \chi_{\{|f| \leq m\}} \right\|_{p(\cdot),q,\alpha} \leq  3 \left\| f-f_{n} \right\|_{p(\cdot),q,\alpha} + 	\left\| f_{n}\chi_{\{|f_{n}| \leq 2m\}} \right\|_{p(\cdot),q,\alpha}.
 \end{align*}
 Since $f_{n}$ is in $V^{0}(L^{p(\cdot)}, \ell^{q})^{\alpha}$, we get 
 \begin{align*}
 	\varlimsup_{m \rightarrow 0} \left\| f \chi_{\{|f| \leq m\}} \right\|_{p(\cdot),q,\alpha} \leq  3 \left\| f-f_{n} \right\|_{p(\cdot),q,\alpha}.
 \end{align*}
 Hence, by letting $n\rightarrow \infty$, we obtain
 \begin{align*}
 	\lim_{m \rightarrow 0} \left\| f \chi_{\{|f| \leq m\}} \right\|_{p(\cdot),q,\alpha} =0.
 \end{align*}
 This shows that $f$ belongs to $V^{0}(L^{p(\cdot)}, \ell^{q})^{\alpha}$. The proof is complete.
 \end{proof}
 
 As a direct consequence of the definitions and the preceding propositions, we obtain the following characterization of the intersection 	$V^{0}(L^{p(\cdot)}, \ell^{q})^{\alpha} \cap V^{\infty}(L^{p(\cdot)}, \ell^{q})^{\alpha}$.
 
 \begin{coro} \label{cor3.4}
 		Let $p(\cdot)$ be a variable exponent on $\mathbb{R}^d$ and $ p(\cdot) \leq  \alpha \leq q \leq \infty$. Then
 	$V^{0}(L^{p(\cdot)}, \ell^{q})^{\alpha} \cap V^{\infty}(L^{p(\cdot)}, \ell^{q})^{\alpha}$ is a closed subset of $(L^{p(\cdot)}, \ell^{q})^{\alpha}$ and we have
 	\begin{align*}
 		V^{0}(L^{p(\cdot)}, \ell^{q})^{\alpha} \cap V^{\infty}(L^{p(\cdot)}, \ell^{q})^{\alpha}=\left\{ f\in (L^{p(\cdot)}, \ell^{q})^{\alpha} : \displaystyle \lim_{m\rightarrow 0} \left\|f- \chi_{\left(m,\frac{1}{m}\right]}(|f|)f\right\|_{p(\cdot),q,\alpha}=0 \right\}.
 	\end{align*}
 \end{coro}
 
 A classical result in variable Lebesgue spaces reads as follows.
 
 \begin{lem}\label{lem3.5} \cite[Proposition 2.18]{Cruz-Fior-2013}.
 	Let $s >0$ and let $p(\cdot) \in \mathcal{P}(\mathbb{R}^d)$ such that $sp(\cdot) \in \mathcal{P}(\mathbb{R}^d)$.  Then, for any $f\in L^{p(\cdot)}$, we have
 	\begin{align*}
 		\||f|^{s}\|_{p(\cdot)}= \|f\|_{sp(\cdot)}^{s}.
 	\end{align*}
 \end{lem}
 
 This result can be extended to variable-exponent Fofana spaces.
 
 \begin{lem}\label{lem3.6} 
 	Let $s >0$ and let $p(\cdot) \in \mathcal{P}(\mathbb{R}^d)$ such that  $sp(\cdot) \in \mathcal{P}(\mathbb{R}^d)$.  Suppose that 
 	$p(\cdot) \leq  \alpha \leq q \leq \infty$. Then, for any $f\in (L^{p(\cdot)}, \ell^{q})^{\alpha}$, we have
 	\begin{align*}
 		\||f|^{s}\|_{p(\cdot),q,\alpha}= \|f\|_{sp(\cdot),sq,s\alpha}^{s}.
 	\end{align*}
 \end{lem}
 
 \begin{proof}
 	We only have to prove the case $q<\infty$, since the case  $q=\infty$ follows more easily. Let  $r$ be a positive real number  and let $k\in \mathbb{Z}^{d}$. Then, by Lemma \ref{lem3.5}, we have
 	\begin{align*}
 			_{r}\||f|^{s}\|_{p(\cdot),q} =\left[\displaystyle \sum_{k\in \mathbb{Z}^{d}} \left\| |f|^{s} \chi_{I_{k}^{r}} \right\|_{p(\cdot)}^{q}\right]^{\frac{1}{q}} =\left[\displaystyle \sum_{k\in \mathbb{Z}^{d}} \left\| f \chi_{I_{k}^{r}} \right\|_{sp(\cdot)}^{sq}\right]^{\frac{1}{sq}\times s} = \, _{r}\|f\|_{sp(\cdot),sq}^{s}.
 	\end{align*}
 Multiplying both sides of the above inequality by $r^{d\left(\frac{1}{\alpha}-\frac{1}{N_{p(\cdot),r}}\right)}$, we get
 	\begin{align*}
 		r^{d\left(\frac{1}{\alpha}-\frac{1}{N_{p(\cdot),r}}\right)}	\, _{r}\||f|^{s}\|_{p(\cdot),q} = \left( 	r^{d\left(\frac{1}{s \alpha}-\frac{1}{s N_{p(\cdot),r}}\right)}	\, _{r}\|f\|_{sp(\cdot),sq}\right)^{s}.
 	\end{align*}
 Therefore, taking the supremum over all $r>0$, we obtain the desired equality.
 \end{proof}

As a consequence of H$\ddot{\mbox{o}}$lder's inequality and Lemma \ref{lem3.6}, we state the following result.
 
 \begin{lem} \label{lem3.7}
 	Let	$\theta \in (0,1)$ and $p_{0}(\cdot), p_{1}(\cdot) \in \mathcal{P}(\mathbb{R}^d)$. Assume that $p_{j}(\cdot) \leq  \alpha_{j}\leq q_{j}\leq\infty$ for all  $j$ in $\{ 0,1\}$, $q_{0}\neq q_{1}$ and
 	\begin{align}\label{eq3.3}
 		\frac{p_{0}(\cdot)}{p_{1}(\cdot)}= \frac{\alpha_{0}}{\alpha_{1}}=\frac{q_{0}}{q_{1}}.
 	\end{align} 
 	Define, for all  $r>0$,
 	\begin{align} \label{eq3.4}
 		\frac{1}{p(\cdot)} = \frac{1-\theta}{p_{0}(\cdot)}+\frac{\theta}{p_{1}(\cdot)}, \frac{1}{N_{p,r}}=\frac{1-\theta}{N_{p_{0},r}}+\frac{\theta}{N_{p_{1},r}}, \frac{1}{\alpha}= \frac{1-\theta}{\alpha_{0}}+\frac{\theta}{\alpha_{1}} \mbox{ and }  \frac{1}{q}=\frac{1-\theta}{q_{0}}+\frac{\theta}{q_{1}}.
 	\end{align}
 	Then, any element $f$ of $(L^{p_{0}(\cdot)},\ell^{q_{0}})^{\alpha_{0}}\cap (L^{p_{1}(\cdot)},\ell^{q_{1}})^{\alpha_{1}}$ belongs to $V^{0}(L^{p(\cdot)}, \ell^{q})^{\alpha} \cap V^{\infty}(L^{p(\cdot)}, \ell^{q})^{\alpha}$ and
 	\begin{align*}
 		\left\|f\right\|_{p(\cdot),q,\alpha} \leq 	\left\|f\right\|_{p_{0}(\cdot),q_{0},\alpha_{0}}^{1-\theta} \left\|f\right\|_{p_{1}(\cdot),q_{1},\alpha_{1}}^{\theta}.
 	\end{align*}
 \end{lem}
 
 \begin{proof}
 		Let $f$ be an element of $(L^{p_{0}(\cdot)},\ell^{q_{0}})^{\alpha_{0}}\cap (L^{p_{1}(\cdot)},\ell^{q_{1}})^{\alpha_{1}}$.\\
 	$\bullet$ Let $r$ be a positive real number. Since 	$\frac{1}{p(\cdot)} = \frac{1-\theta}{p_{0}(\cdot)}+\frac{\theta}{p_{1}(\cdot)}$, H$\ddot{\mbox{o}}$lder's inequality for variable exponent (see \cite[Corollary 2.28]{Cruz-Fior-2013}) implies that there exists a constant $C$ such that, for all  $k\in \mathbb{Z}^{d}$,
 	\begin{align*}
 		\|f\chi_{I_{k}^{r}}\|_{p(\cdot)} \leq C \;\|f\chi_{I_{k}^{r}}\|_{p_{0}(\cdot)}^{1-\theta} \|f\chi_{I_{k}^{r}}\|_{p_{1}(\cdot)}^{\theta}.
 	\end{align*}
 	Therefore, since $\frac{1}{q}=\frac{1-\theta}{q_{0}}+\frac{\theta}{q_{1}}$, if $q<\infty$ then we may apply H$\ddot{\mbox{o}}$lder's inequality to get
 	\begin{align*}
 		_{r}\|f\|_{p(\cdot),q} \leq C \left[\sum_{k\in \mathbb{Z}^{d}} 	\|f\chi_{I_{k}^{r}}\|_{p_{0}(\cdot)}^{q_{0}} \right]^{\frac{1-\theta}{q_{0}}} \left[\sum_{k\in \mathbb{Z}^{d}} 	\|f\chi_{I_{k}^{r}}\|_{p_{1}(\cdot)}^{q_{1}} \right]^{\frac{\theta}{q_{1}}} = C \; _{r}\|f\|_{p_{0}(\cdot),q_{0}}^{1-\theta} \, _{r}\|f\|_{p_{1}(\cdot),q_{1}}^{\theta}.
 	\end{align*}
 	Likwise, if $q=\infty$ then
 	\begin{align*}
 		_{r}\|f\|_{p(\cdot),\infty}= \sup_{k\in \mathbb{Z}^{d}} \|f\chi_{I_{k}^{r}}\|_{p(\cdot)} \leq  \sup_{k\in \mathbb{Z}^{d}} C \;\|f\chi_{I_{k}^{r}}\|_{p_{0}(\cdot)}^{1-\theta} \|f\chi_{I_{k}^{r}}\|_{p_{1}(\cdot)}^{\theta}= C\, 	_{r}\|f\|_{p_{0}(\cdot),\infty}^{1-\theta} \,	_{r}\|f\|_{p_{1}(\cdot),\infty}^{\theta}.
 	\end{align*}
 	Thus, we get
 		\begin{align*}
 		r^{d\left(\frac{1}{\alpha}-\frac{1}{N_{p,r}}\right)}\,		_{r}\|f\|_{p(\cdot),q}
 		\leq C \left(  r^{d\left(\frac{1}{\alpha_{0}}-\frac{1}{N_{p_{0},r}}\right)} \; _{r}\|f\|_{p_{0}(\cdot),q_{0}}\right)^{1-\theta} \left( r^{d\left(\frac{1}{\alpha_{1}}-\frac{1}{N_{p_{1},r}}\right)}\; _{r}\|f\|_{p_{1}(\cdot),q_{1}} \right)^{\theta}.
 	\end{align*}
 	Therefore, taking the supremum over all $r>0$, we obtain 
 	\begin{align*}
 		\left\|f\right\|_{p(\cdot),q,\alpha} \leq 	\left\|f\right\|_{p_{0}(\cdot),q_{0},\alpha_{0}}^{1-\theta} \left\|f\right\|_{p_{1}(\cdot),q_{1},\alpha_{1}}^{\theta}.
 	\end{align*}
 	$\bullet$ Without loss of generality, we may assume that $p_{0}(\cdot)< p_{1}(\cdot)$ and so $p_{0}(\cdot)<p(\cdot)< p_{1}(\cdot)$. Let $m$ be a real number such that $0<m<1$. It is easy to see that
 	\begin{align*}
 		f- \chi_{\left(m,\frac{1}{m}\right]}(|f|)f= f\chi_{\left\{|f|\leq m \right\}\cup \left\{|f|> \frac{1}{m}\right\} }
 	\end{align*}
 	and, for all $x\in \mathbb{R}^{d}$,
 	\begin{align*}
 		\left\{
 		\begin{array}{rl}
 			& |f(x)|\leq m  \Longrightarrow \left(\frac{1}{m} |f(x)|\right)^{p(\cdot)} \leq \left(\frac{1}{m} |f(x)|\right)^{p_{0}(\cdot)} \Longrightarrow|f(x)|\leq m^{1-\frac{p_{0}(\cdot)}{p(\cdot)}} |f(x)|^{\frac{p_{0}(\cdot)}{p(\cdot)}}, \\\\
 			& |f(x)|> \frac{1}{m}  \Longrightarrow \left(m |f(x)|\right)^{p_{1}(\cdot)} > \left(m |f(x)|\right)^{p(\cdot)} \Longrightarrow |f(x)|< m^{\frac{p_{1}(\cdot)}{p(\cdot)}-1} |f(x)|^{\frac{p_{1}(\cdot)}{p(\cdot)}}.
 		\end{array}
 		\right. 
 	\end{align*}
 	Hence, we obtain
 	\begin{align}\label{eq3.5'}
 		\left| f- \chi_{\left(m,\frac{1}{m}\right]}(|f|)f \right| \leq  m^{1-\frac{p_{0}(\cdot)}{p(\cdot)}} |f|^{\frac{p_{0}(\cdot)}{p(\cdot)}}+ m^{\frac{p_{1}(\cdot)}{p(\cdot)}-1} |f|^{\frac{p_{1}(\cdot)}{p(\cdot)}}.
 	\end{align}
 		Observe that, the equality \eqref{eq3.3} implies that, for every $j\in \left\{0,1\right\}$,
 		\begin{align*}
 			\frac{p_{j}(\cdot)}{p(\cdot)}= \frac{\alpha_{j}}{\alpha}=\frac{q_{j}}{q}.
 		\end{align*}
 Therefore, applying Minkowski's inequality to \eqref{eq3.5'} and  using Lemma \ref{lem3.6}, we get
 	\begin{align*}
 		\left\|	f- \chi_{\left(m,\frac{1}{m}\right]}(|f|)f \right\|_{p(\cdot),q,\alpha} \leq  m^{1-\frac{p_{0}(\cdot)}{p(\cdot)}}  \,\| f \|_{p_{0}(\cdot),q_{0},\alpha_{0}}^{\frac{p_{0}(\cdot)}{p(\cdot)}} +   m^{\frac{p_{1}(\cdot)}{p(\cdot)}-1} \, \| f \|_{p_{1}(\cdot),q_{1},\alpha_{1}}^{\frac{p_{1}(\cdot)}{p(\cdot)}}.
 	\end{align*}
 	Moreover, since $p_{0}(\cdot)<p(\cdot)< p_{1}(\cdot)$, we obtain
 	\begin{align*}
 		\displaystyle \lim_{m\rightarrow 0} \left\|f- \chi_{\left(m,\frac{1}{m}\right]}(|f|)f\right\||_{p(\cdot),q,\alpha}=0.
 	\end{align*}
 Hence, Corollary \ref{cor3.4} asserts that $f$ belongs to $V^{0}(L^{p(\cdot)}, \ell^{q})^{\alpha} \cap V^{\infty}(L^{p(\cdot)}, \ell^{q})^{\alpha}$.  The proof is complete.
 \end{proof}


\section{Main results} \label{Sec4}

We begin this section with a result that identifies Calder\'{o}n's product of variable-exponent Fofana spaces with the corresponding intermediate spaces.

\begin{prop} \label{P4.1}
	Let	$\theta \in (0,1)$ and $p_{0}(\cdot), p_{1}(\cdot) \in \mathcal{P}(\mathbb{R}^d)$. Assume that $p_{j}(\cdot) \leq  \alpha_{j}\leq q_{j}\leq\infty$ for all  $j$ in $\{ 0,1\}$, $q_{0}\neq q_{1}$ and
	\begin{align*}
		\frac{p_{0}(\cdot)}{p_{1}(\cdot)}= \frac{\alpha_{0}}{\alpha_{1}}=\frac{q_{0}}{q_{1}}.
	\end{align*} 
	Define, for all  $r>0$,
	\begin{align*} 
		\frac{1}{p(\cdot)} = \frac{1-\theta}{p_{0}(\cdot)}+\frac{\theta}{p_{1}(\cdot)}, \frac{1}{N_{p,r}}=\frac{1-\theta}{N_{p_{0},r}}+\frac{\theta}{N_{p_{1},r}}, \frac{1}{\alpha}= \frac{1-\theta}{\alpha_{0}}+\frac{\theta}{\alpha_{1}} \mbox{ and }  \frac{1}{q}=\frac{1-\theta}{q_{0}}+\frac{\theta}{q_{1}}.
	\end{align*}
	 Then, we have
	\begin{align*}
		\left[(L^{p_{0}(\cdot)},\ell^{q_{0}})^{\alpha_{0}}\right]^{1-\theta} \left[(L^{p_{1}(\cdot)},\ell^{q_{1}})^{\alpha_{1}} \right]^{\theta} = (L^{p(\cdot)}, \ell^{q})^{\alpha}.
	\end{align*}
\end{prop}

\begin{proof}
	Let $f$ be an element of $\left[(L^{p_{0}(\cdot)},\ell^{q_{0}})^{\alpha_{0}}\right]^{1-\theta} \left[(L^{p_{1}(\cdot)},\ell^{q_{1}})^{\alpha_{1}} \right]^{\theta}$. Then, there exists $f_{0}\in (L^{p_{0}(\cdot)},\ell^{q_{0}})^{\alpha_{0}} $ and $f_{1}\in (L^{p_{1}(\cdot)},\ell^{q_{1}})^{\alpha_{1}}$ such that $|f|\leq |f_{0}|^{1-\theta}|f_{1}|^{\theta}$. Thus, H$\ddot{\mbox{o}}$lder's inequality for variable exponents implies that there exists a constant $C$ such that, for all real numbers $r>0$ and for all $k\in \mathbb{Z}^{d}$, we have
	\begin{align*}
		r^{d\left(\frac{1}{\alpha}-\frac{1}{N_{p,r}}\right)}	\|f\chi_{I_{k}^{r}}\|_{p(\cdot)}
		\leq C \left(  r^{d\left(\frac{1}{\alpha_{0}}-\frac{1}{N_{p_{0},r}}\right)} \|f_{0}\chi_{I_{k}^{r}}\|_{p_{0}(\cdot)}\right)^{1-\theta} \left( r^{d\left(\frac{1}{\alpha_{1}}-\frac{1}{N_{p_{1},r}}\right)} \|f_{1}\chi_{I_{k}^{r}}\|_{p_{1}(\cdot)}\right)^{\theta}.
	\end{align*}
	Therefore, summing over $k$, using H$\ddot{\mbox{o}}$lder's inequality, and then taking the supremum over $r>0$, we obtain
	\begin{align*}
		\left\|f\right\|_{p(\cdot),q,\alpha} \leq C	\left\|f_{0}\right\|_{p_{0}(\cdot),q_{0},\alpha_{0}}^{1-\theta} \left\|f_{1}\right\|_{p_{1}(\cdot),q_{1},\alpha_{1}}^{\theta} < \infty.
	\end{align*}
	This shows that $f\in (L^{p(\cdot)}, \ell^{q})^{\alpha}$ and so
	\begin{align} \label{eq3.5}
		\left[(L^{p_{0}(\cdot)},\ell^{q_{0}})^{\alpha_{0}}\right]^{1-\theta} \left[(L^{p_{1}(\cdot)},\ell^{q_{1}})^{\alpha_{1}} \right]^{\theta} \subset (L^{p(\cdot)}, \ell^{q})^{\alpha}.
	\end{align}
	Conversely, let $f$ be an element of $(L^{p(\cdot)}, \ell^{q})^{\alpha}$. By virtue of Lemma \ref{lem3.6}, for all $j\in \{0, 1\}$, 
	\begin{align*}
		\left\||f|^{\frac{p(\cdot{})}{p_{j}(\cdot{})}}\right\|_{p_{j}(\cdot),q_{j},\alpha_{j}}= \left\|f\right\|_{p(\cdot),q,\alpha}^{\frac{p(\cdot{})}{p_{j}(\cdot{})}}
	\end{align*}
	and subsequently,	$f_{j}=|f|^{\frac{p(\cdot{})}{p_{j}(\cdot{})}}$ belongs to $(L^{p_{j}(\cdot)},\ell^{q_{j}})^{\alpha_{j}}$. Furthermore, observe that $|f_{0}|^{1-\theta} |f_{1}|^{\theta} =|f|$. Hence, $f$ belongs to $\left[(L^{p_{0}(\cdot)},\ell^{q_{0}})^{\alpha_{0}}\right]^{1-\theta} \left[(L^{p_{1}(\cdot)},\ell^{q_{1}})^{\alpha_{1}} \right]^{\theta}$ and so
	\begin{align}\label{eq3.6}
		(L^{p(\cdot)}, \ell^{q})^{\alpha} \subset \left[(L^{p_{0}(\cdot)},\ell^{q_{0}})^{\alpha_{0}}\right]^{1-\theta} \left[(L^{p_{1}(\cdot)},\ell^{q_{1}})^{\alpha_{1}} \right]^{\theta}.
	\end{align}
	Combining \eqref{eq3.5} and  \eqref{eq3.6}, we obtain the desired result.
\end{proof}

We are now in a position to state and prove interpolation results of this paper.

\begin{theo} \label{Th4.2}
Let	$\theta \in (0,1)$ and $p_{0}(\cdot), p_{1}(\cdot) \in \mathcal{P}(\mathbb{R}^d)$. Assume that $p_{j}(\cdot) \leq  \alpha_{j}\leq q_{j}\leq\infty$ for all  $j$ in $\{ 0,1\}$, $q_{0}\neq q_{1}$ and
\begin{align*}
	\frac{p_{0}(\cdot)}{p_{1}(\cdot)}= \frac{\alpha_{0}}{\alpha_{1}}=\frac{q_{0}}{q_{1}}.
\end{align*} 
Define, for all $r>0$,
	\begin{align*} 
	\frac{1}{p(\cdot)} = \frac{1-\theta}{p_{0}(\cdot)}+\frac{\theta}{p_{1}(\cdot)}, \frac{1}{N_{p,r}}=\frac{1-\theta}{N_{p_{0},r}}+\frac{\theta}{N_{p_{1},r}}, \frac{1}{\alpha}= \frac{1-\theta}{\alpha_{0}}+\frac{\theta}{\alpha_{1}} \mbox{ and }  \frac{1}{q}=\frac{1-\theta}{q_{0}}+\frac{\theta}{q_{1}}.
\end{align*}
Then, we have
\begin{align*}
	\left[(L^{p_{0}(\cdot)},\ell^{q_{0}})^{\alpha_{0}},(L^{p_{1}(\cdot)},\ell^{q_{1}})^{\alpha_{1}}\right]_{\theta}= V^{0}(L^{p(\cdot)}, \ell^{q})^{\alpha} \cap V^{\infty}(L^{p(\cdot)}, \ell^{q})^{\alpha}.
\end{align*}
\end{theo}

\begin{proof}
	 Since $(L^{p_{j}(\cdot)},\ell^{q_{j}})^{\alpha_{j}}\,(j=0,1)$ are Banach function spaces, \eqref{eq 2.1} and Proposition \ref{P4.1} imply that
	\begin{align}\label{eq4.3}
		\left[(L^{p_{0}(\cdot)},\ell^{q_{0}})^{\alpha_{0}},(L^{p_{1}(\cdot)},\ell^{q_{1}})^{\alpha_{1}}\right]_{\theta}= \overline{(L^{p_{0}(\cdot)},\ell^{q_{0}})^{\alpha_{0}} \cap (L^{p_{1}(\cdot)},\ell^{q_{1}})^{\alpha_{1}}}^{(L^{p(\cdot)}, \ell^{q})^{\alpha}}.
	\end{align}
	Furthermore, Lemma \ref{lem3.7} asserts that $(L^{p_{0}(\cdot)},\ell^{q_{0}})^{\alpha_{0}} \cap (L^{p_{1}(\cdot)},\ell^{q_{1}})^{\alpha_{1}}$ is included in  $V^{0}(L^{p(\cdot)}, \ell^{q})^{\alpha} \cap V^{\infty}(L^{p(\cdot)}, \ell^{q})^{\alpha}$, which is a closed subspace of  $(L^{p(\cdot)}, \ell^{q})^{\alpha}$ (see Corollary \ref{cor3.4}). Thus
	\begin{align*}
	\overline{(L^{p_{0}(\cdot)},\ell^{q_{0}})^{\alpha_{0}} \cap (L^{p_{1}(\cdot)},\ell^{q_{1}})^{\alpha_{1}}}^{(L^{p(\cdot)}, \ell^{q})^{\alpha}} \subset V^{0}(L^{p(\cdot)}, \ell^{q})^{\alpha} \cap V^{\infty}(L^{p(\cdot)}, \ell^{q})^{\alpha}
	\end{align*}
	and therefore, by \eqref{eq4.3}, we obtain
	\begin{align}\label{eq4.4}
			\left[(L^{p_{0}(\cdot)},\ell^{q_{0}})^{\alpha_{0}},(L^{p_{1}(\cdot)},\ell^{q_{1}})^{\alpha_{1}}\right]_{\theta} \subset V^{0}(L^{p(\cdot)}, \ell^{q})^{\alpha} \cap V^{\infty}(L^{p(\cdot)}, \ell^{q})^{\alpha}.
	\end{align}

Without loss of generality, we may assume that $p_{0}(\cdot)< p_{1}(\cdot)$ and so $p_{0}(\cdot)<p(\cdot)< p_{1}(\cdot)$.
	Let $f$ be an element of $V^{0}(L^{p(\cdot)}, \ell^{q})^{\alpha} \cap V^{\infty}(L^{p(\cdot)}, \ell^{q})^{\alpha}$
	and consider a real number $\epsilon>0$.
	Since $f \in V^{0}(L^{p(\cdot)}, \ell^{q})^{\alpha}$, there exists a real number $m_{\epsilon}>0$ such that for all integers $m\leq m_{\epsilon}$, we have 
	\begin{align}\label{eq4.5}
		\left\|f\chi_{\{|f|\leq m\}}\right\|_{p(\cdot),q,\alpha}<\frac{\epsilon}{2}.
	\end{align}
	Likwise, since $f \in  V^{\infty}(L^{p(\cdot)}, \ell^{q})^{\alpha}$, there exists $g_{\epsilon}$ in $L^{\infty}\cap(L^{q},\ell^{p})^{\alpha}$ such that
	\begin{align}\label{eq4.6}
		\left\|f-g_{\epsilon} \right\|_{q,p,\alpha}<\frac{\epsilon}{2}.
	\end{align}
	Let us consider a real number $m>0$. Then, since $p_{0}(\cdot)<p(\cdot)$, we have
	\begin{align*}
		\left\|f\chi_{\{|f|>m\}}\right\|_{p_{0}(\cdot),q_{0},\alpha_{0}} &= \left\| |f|^{1-\frac{p(\cdot)}{p_{0}(\cdot)}} |f|^{\frac{p(\cdot)}{p_{0}(\cdot)}}\chi_{\{|f|>m\}}\right\|_{p_{0}(\cdot),q_{0},\alpha_{0}} \\
		&\leq m^{1-\frac{p(\cdot)}{p_{0}(\cdot)}} \left\| |f|^{\frac{p(\cdot)}{p_{0}(\cdot)}} \right\|_{p_{0}(\cdot),q_{0},\alpha_{0}}.
	\end{align*}
	Therefore, by virtue of Lemma \ref{lem3.6}, we get 
\begin{align}\label{eq4.7}
	\left\|f\chi_{\{|f|>m\}}\right\|_{p_{0}(\cdot),q_{0},\alpha_{0}} 
	\leq m^{1-\frac{p(\cdot)}{p_{0}(\cdot)}} \left\| f \right\|_{p(\cdot),q,\alpha}^{\frac{p(\cdot)}{p_{0}(\cdot)}} <\infty.
\end{align}
	Furthermore, since $p(\cdot)<p_{1}(\cdot)$, another application of Lemma \ref{lem3.6} leads to
	\begin{align}\label{eq4.8}
		\left\|g_{\epsilon}\right\|_{p_{1}(\cdot),q_{1},\alpha_{1}} 
		\leq  \|g_{\epsilon}\|_{\infty}^{1-\frac{p(\cdot)}{p_{1}(\cdot)}}  \left\|  g_{\epsilon}\right\|_{p_{1}(\cdot),q_{1},\alpha_{1}}^{\frac{p(\cdot)}{p_{1}(\cdot)}} <\infty.
	\end{align}

	Let us set, for all $m>0$,
	\begin{align*}
		g_{\epsilon,m}=
		\left\{
		\begin{array}{rl}
			& f\chi_{\{|f|> m\}}\;\;\;\;\mbox{ if } \;\;|g_{\epsilon}|>|f|,\\
			&g_{\epsilon}\chi_{\{|f|> m\}}\;\;\;\mbox{ if } \;\;|g_{\epsilon}|\leq |f|.
		\end{array}
		\right.
	\end{align*}
	It is easy to see that
	$$|g_{\epsilon,m}|=\min\left( |f|\chi_{\{|f|> m\}}, \;|g_{\epsilon}|\chi_{\{|f|> m\}} \right)\leq  \min\left( |f|\chi_{\{|f|> m\}}, \;|g_{\epsilon}|\right).$$
	Thus, in accordance with \eqref{eq4.7} and \eqref{eq4.8}, $g_{\epsilon,m}$ belongs to  $(L^{p_{0}(\cdot)},\ell^{q_{0}})^{\alpha_{0}} \cap (L^{p_{1}(\cdot)},\ell^{q_{1}})^{\alpha_{1}}$.\\
	Moreover, we have
	\begin{align*}
		\left| f-g_{\epsilon,m}\right|
		&\leq \left| f-f\chi_{\{|f|> m\}}\right|+ \left|f\chi_{\{|f|> m\}}-g_{\epsilon,m}\right|\\
		&\leq \left| f\chi_{\{|f|\leq m\}}\right|+ \left|f-g_{\epsilon}\right|.
	\end{align*}
	Consequently, for all $m\leq m_{\epsilon}$, \eqref{eq4.5} and \eqref{eq4.6} lead to
	\begin{align*}
		\left\|f-g_{\epsilon,m} \right\|_{p(\cdot),q,\alpha} \leq  \left\|f\chi_{\{|f|\leq m\}} \right\|_{p(\cdot),q,\alpha} + \left\|f-g_{\epsilon} \right\|_{p(\cdot),q,\alpha}<\epsilon.
	\end{align*}
	Hence, $f$ belongs to the closure  of $(L^{p_{0}(\cdot)},\ell^{q_{0}})^{\alpha_{0}} \cap (L^{p_{1}(\cdot)},\ell^{q_{1}})^{\alpha_{1}}$ in $(L^{p(\cdot)},\ell^{q})^{\alpha}$. This shows that
	\begin{align*}
 V^{0}(L^{p(\cdot)}, \ell^{q})^{\alpha} \cap V^{\infty}(L^{p(\cdot)}, \ell^{q})^{\alpha} \subset 	\overline{(L^{p_{0}(\cdot)},\ell^{q_{0}})^{\alpha_{0}} \cap (L^{p_{1}(\cdot)},\ell^{q_{1}})^{\alpha_{1}}}^{(L^{p(\cdot)}, \ell^{q})^{\alpha}}.
	\end{align*}
	Therefore, by \eqref{eq4.3}, we obtain 
	\begin{align}\label{eq4.9}
		 V^{0}(L^{p(\cdot)}, \ell^{q})^{\alpha} \cap V^{\infty}(L^{p(\cdot)}, \ell^{q})^{\alpha} \subset	\left[(L^{p_{0}(\cdot)},\ell^{q_{0}})^{\alpha_{0}},(L^{p_{1}(\cdot)},\ell^{q_{1}})^{\alpha_{1}}\right]_{\theta}.
	\end{align}
	The conjunction of \eqref{eq4.4} and \eqref{eq4.9} provides the desired result.
\end{proof}

We obtain the following result.

\begin{theo}  \label{Th4.3}
	Keep the same assumptions as in Theorem \ref{Th4.2}. Then, we have
	\begin{align*}
		\left[(L^{p_{0}(\cdot)},\ell^{q_{0}})^{\alpha_{0}},(L^{p_{1}(\cdot)},\ell^{q_{1}})^{\alpha_{1}}\right]^{\theta}= (L^{p(\cdot)}, \ell^{q})^{\alpha}.
	\end{align*}
\end{theo}

\begin{proof}
	Let $(f_{n})_{n\geq 1}$ be a Cauchy sequence in the unit ball of  the Calder\'{o}n product $\left[(L^{p_{0}(\cdot)},\ell^{q_{0}})^{\alpha_{0}}\right]^{1-\theta} \left[(L^{p_{1}(\cdot)},\ell^{q_{1}})^{\alpha_{1}} \right]^{\theta}$. Then, by virtue of Proposition \ref{P4.1}, $(f_{n})_{n\geq 1}$ is also a Cauchy sequence in $(L^{p(\cdot)}, \ell^{q})^{\alpha}$ and we have
	$\left\|f_{n} \right\|_{p(\cdot),q,\alpha}<1$, for all  integers $n\geq 1$. Since $(L^{p(\cdot)}, \ell^{q})^{\alpha}$ is continuously included in $(L^{p_{0}(\cdot)},\ell^{q_{0}})^{\alpha_{0}}+(L^{p_{1}(\cdot)},\ell^{q_{1}})^{\alpha_{1}}$, $(f_{n})_{n\geq 1}$ is a Cauchy sequence in $(L^{p_{0}(\cdot)},\ell^{q_{0}})^{\alpha_{0}}+(L^{p_{1}(\cdot)},\ell^{q_{1}})^{\alpha_{1}}$. Thus, $(f_{n})_{n\geq 1}$ converges to a function $f$ in $(L^{p_{0}(\cdot)},\ell^{q_{0}})^{\alpha_{0}}+(L^{p_{1}(\cdot)},\ell^{q_{1}})^{\alpha_{1}}$. Consequently, there exists a subsequence $(f_{n_{k}})_{k\geq 1}$ of $(f_{n})_{n\geq 1}$ such that $\displaystyle \lim_{k\rightarrow \infty}f_{n_{k}}(x)=f(x)$, for almost every $x\in \mathbb{R}^{d}$. Hence, by the Fatou lemma, we have
	\begin{align*}
		\|f\|_{\left[(L^{p_{0}(\cdot)},\ell^{q_{0}})^{\alpha_{0}}\right]^{1-\theta} \left[(L^{p_{1}(\cdot)},\ell^{q_{1}})^{\alpha_{1}} \right]^{\theta}} = \|f \|_{p(\cdot),q,\alpha} \leq \displaystyle \liminf_{k\rightarrow \infty} \|f_{n_{k}}\|_{p(\cdot),q,\alpha}  < 1.
	\end{align*}
	This proves that $f$ is in the unit ball of $\left[(L^{p_{0}(\cdot)},\ell^{q_{0}})^{\alpha_{0}}\right]^{1-\theta} \left[(L^{p_{1}(\cdot)},\ell^{q_{1}})^{\alpha_{1}} \right]^{\theta}$. Thus, the unit ball of $\left[(L^{p_{0}(\cdot)},\ell^{q_{0}})^{\alpha_{0}}\right]^{1-\theta} \left[(L^{p_{1}(\cdot)},\ell^{q_{1}})^{\alpha_{1}} \right]^{\theta}$ is closed in  $(L^{p_{0}(\cdot)},\ell^{q_{0}})^{\alpha_{0}}+(L^{p_{1}(\cdot)},\ell^{q_{1}})^{\alpha_{1}}$. Then, in accordance with \eqref{eq 2.3} and Proposition \ref{P4.1}, we obtain the desired result.
\end{proof}

The following result holds.

\begin{theo} \label{Th4.4}
	Keep the same assumptions as in Theorem \ref{Th4.2}. Then, we have
	\begin{align*}
		\left[(L^{p_{0}(\cdot)},\ell^{q_{0}})^{\alpha_{0}}_{\rm a},(L^{p_{1}(\cdot)},\ell^{q_{1}})^{\alpha_{1}}_{\rm a}\right]_{\theta}=\left[(L^{p_{0}(\cdot)},\ell^{q_{0}})^{\alpha_{0}}_{\rm a},(L^{p_{1}(\cdot)},\ell^{q_{1}})^{\alpha_{1}}\right]_{\theta}= (L^{p(\cdot)}, \ell^{q})^{\alpha}_{\rm a}.
	\end{align*}
\end{theo}

\begin{proof}
	Since, for each $j\in \{0, 1\}$, the space $(L^{p_{j}(\cdot)},\ell^{q_{j}})^{\alpha_{j}}_{\rm a}$ has an absolutely continuous norm, Remark \ref{rem 2.6} implies that 
	\begin{align} \label{eq4.10}
		\left[(L^{p_{0}(\cdot)},\ell^{q_{0}})^{\alpha_{0}}_{\rm a},(L^{p_{1}(\cdot)},\ell^{q_{1}})^{\alpha_{1}}_{\rm a}\right]_{\theta} = \left[(L^{p_{0}(\cdot)},\ell^{q_{0}})^{\alpha_{0}}_{\rm a}\right]^{1-\theta} \left[(L^{p_{1}(\cdot)},\ell^{q_{1}})^{\alpha_{1}}_{\rm a} \right]^{\theta}
	\end{align}
	and
	\begin{align} \label{eq4.11}
		\left[(L^{p_{0}(\cdot)},\ell^{q_{0}})^{\alpha_{0}}_{\rm a},(L^{p_{1}(\cdot)},\ell^{q_{1}})^{\alpha_{1}}\right]_{\theta} = \left[(L^{p_{0}(\cdot)},\ell^{q_{0}})^{\alpha_{0}}_{\rm a}\right]^{1-\theta} \left[(L^{p_{1}(\cdot)},\ell^{q_{1}})^{\alpha_{1}} \right]^{\theta}.
	\end{align}
	
	Moreover, it follows from Remark \ref{rem 2.6} that the spaces $\left[(L^{p_{0}(\cdot)},\ell^{q_{0}})^{\alpha_{0}}_{\rm a}\right]^{1-\theta} \left[(L^{p_{1}(\cdot)},\ell^{q_{1}})^{\alpha_{1}}_{\rm a} \right]^{\theta}$ and $\left[(L^{p_{0}(\cdot)},\ell^{q_{0}})^{\alpha_{0}}_{\rm a}\right]^{1-\theta} \left[(L^{p_{1}(\cdot)},\ell^{q_{1}})^{\alpha_{1}} \right]^{\theta}$  have absolutely continuous norms, and therefore    Proposition \ref{P4.1} implies that these two spaces are embedded in $(L^{p(\cdot)}, \ell^{q})^{\alpha}_{\rm a}$.
	
	Let $f$ be an element of $(L^{p(\cdot)}, \ell^{q})^{\alpha}_{\rm a}$ and let $f_{j}=|f|^{\frac{p(\cdot)}{p_{j}(\cdot)}}$, for any $j\in \{0, 1\}$. In accordance with Lemma \ref{lem3.6}, for each $j\in \{0, 1\}$, the function $f_{j}$ belongs to $(L^{p_{j}(\cdot)},\ell^{q_{j}})^{\alpha_{j}}_{\rm a}$. Meanwhile, observe that $|f|=|f_{0}|^{1-\theta}|f_{1}|^{\theta}$. Consequently, $f$ belongs to the spaces $\left[(L^{p_{0}(\cdot)},\ell^{q_{0}})^{\alpha_{0}}_{\rm a}\right]^{1-\theta} \left[(L^{p_{1}(\cdot)},\ell^{q_{1}})^{\alpha_{1}}_{\rm a} \right]^{\theta}$ and $\left[(L^{p_{0}(\cdot)},\ell^{q_{0}})^{\alpha_{0}}_{\rm a}\right]^{1-\theta} \left[(L^{p_{1}(\cdot)},\ell^{q_{1}})^{\alpha_{1}} \right]^{\theta}$. This shows that the space $(L^{p(\cdot)}, \ell^{q})^{\alpha}_{\rm a}$ is included in both $\left[(L^{p_{0}(\cdot)},\ell^{q_{0}})^{\alpha_{0}}_{\rm a}\right]^{1-\theta} \left[(L^{p_{1}(\cdot)},\ell^{q_{1}})^{\alpha_{1}}_{\rm a} \right]^{\theta}$ and $\left[(L^{p_{0}(\cdot)},\ell^{q_{0}})^{\alpha_{0}}_{\rm a}\right]^{1-\theta} \left[(L^{p_{1}(\cdot)},\ell^{q_{1}})^{\alpha_{1}} \right]^{\theta}$.
	
	In summary, we have proved that 
	\begin{align} \label{eq4.12}
		\left[(L^{p_{0}(\cdot)},\ell^{q_{0}})^{\alpha_{0}}_{\rm a}\right]^{1-\theta} \left[(L^{p_{1}(\cdot)},\ell^{q_{1}})^{\alpha_{1}}_{\rm a} \right]^{\theta} = \left[(L^{p_{0}(\cdot)},\ell^{q_{0}})^{\alpha_{0}}_{\rm a}\right]^{1-\theta} \left[(L^{p_{1}(\cdot)},\ell^{q_{1}})^{\alpha_{1}} \right]^{\theta} = (L^{p(\cdot)}, \ell^{q})^{\alpha}_{\rm a}.
	\end{align}
	The conjunction of \eqref{eq4.10}, \eqref{eq4.11} and \eqref{eq4.12} leads to the desired result.
\end{proof}

Now, we state the following interpolation result concerning preduals of variable-exponent Fofana spaces.

\begin{theo}  \label{Th4.5}
Keep the same assumptions as in Theorem \ref{Th4.2}. Moreover, suppose that $p_{j}(\cdot) \neq  \alpha_{j}\neq q_{j}$ for each  $j$ in $\{ 0,1\}$. Then, we have
	\begin{align*}
	\left[\mathcal{H}\left(p_{0}^{\prime}(\cdot),q_{0}^{\prime},\alpha_{0}^{\prime}\right), \mathcal{H}\left(p_{1}^{\prime}(\cdot),q_{1}^{\prime},\alpha_{1}^{\prime}\right)\right]_{\theta}=		\left[\mathcal{H}\left(p_{0}^{\prime}(\cdot),q_{0}^{\prime},\alpha_{0}^{\prime}\right), \mathcal{H}\left(p_{1}^{\prime}(\cdot),q_{1}^{\prime},\alpha_{1}^{\prime}\right) \right]^{\theta}= \mathcal{H}\left(p^{\prime}(\cdot),q^{\prime},\alpha^{\prime}\right).
\end{align*}	
\end{theo}

\begin{proof}
Since,  for each  $j$ in $\{ 0,1\}$, the space $\mathcal{H}\left(p_{0}^{\prime}(\cdot),q_{0}^{\prime},\alpha_{0}^{\prime}\right)$ has an absolutely continuous norm (see Proposition \ref{P2.9}), \eqref{eq 2.2} implies that
\begin{align*}
		\left[\mathcal{H}\left(p_{0}^{\prime}(\cdot),q_{0}^{\prime},\alpha_{0}^{\prime}\right), \mathcal{H}\left(p_{1}^{\prime}(\cdot),q_{1}^{\prime},\alpha_{1}^{\prime}\right)\right]_{\theta} =	\mathcal{H}\left(p_{0}^{\prime}(\cdot),q_{0}^{\prime},\alpha_{0}^{\prime}\right)^{1-\theta} \mathcal{H}\left(p_{1}^{\prime}(\cdot),q_{1}^{\prime},\alpha_{1}^{\prime}\right)^{\theta}.
\end{align*}
Therefore, applying \eqref{eq2.5}, Proposition \ref{prop 2.7}, and Proposition \ref{P4.1} successively, we obtain
\begin{align*}
	\left[\mathcal{H}\left(p_{0}^{\prime}(\cdot),q_{0}^{\prime},\alpha_{0}^{\prime}\right), \mathcal{H}\left(p_{1}^{\prime}(\cdot),q_{1}^{\prime},\alpha_{1}^{\prime}\right)\right]_{\theta} &=	\left(\left[(L^{p_{0}(\cdot)},\ell^{q_{0}})^{\alpha_{0}}\right]^{\prime}\right)^{1-\theta} \left(\left[(L^{p_{1}(\cdot)},\ell^{q_{1}})^{\alpha_{1}}\right]^{\prime}\right)^{\theta}\\
	&= \left(\left[(L^{p_{0}(\cdot)},\ell^{q_{0}})^{\alpha_{0}}\right]^{1-\theta}
	\left[(L^{p_{1}(\cdot)},\ell^{q_{1}})^{\alpha_{1}}\right]^{\theta}\right)^{\prime}\\
	&= \left[(L^{p(\cdot)},\ell^{q})^{\alpha}\right]^{\prime}.
\end{align*}
Then another application of \eqref{eq2.5} leads to 
\begin{align} \label{eq4.13}
	\left[\mathcal{H}\left(p_{0}^{\prime}(\cdot),q_{0}^{\prime},\alpha_{0}^{\prime}\right), \mathcal{H}\left(p_{1}^{\prime}(\cdot),q_{1}^{\prime},\alpha_{1}^{\prime}\right)\right]_{\theta}=		 \mathcal{H}\left(p^{\prime}(\cdot),q^{\prime},\alpha^{\prime}\right).
\end{align}	

Furthermore, in accordance with \eqref{eq2.6}, we have
\begin{align} \label{eq4.14}
		\left[\mathcal{H}\left(p_{0}^{\prime}(\cdot),q_{0}^{\prime},\alpha_{0}^{\prime}\right), \mathcal{H}\left(p_{1}^{\prime}(\cdot),q_{1}^{\prime},\alpha_{1}^{\prime}\right)\right]^{\theta} = 	\left[ \left((L^{p_{0}(\cdot)},\ell^{q_{0}})^{\alpha_{0}}_{\rm a}\right)^{*}, \left((L^{p_{1}(\cdot)},\ell^{q_{1}})^{\alpha_{1}}_{\rm a}\right)^{*} \right]^{\theta}.
\end{align}
Since $\overline{\mathcal{C}_{\rm c}^{\infty} }= (L^{p_{j}(\cdot)},\ell^{q_{j}})^{\alpha_{j}}_{\rm a}$ for all  $j$ in $\{ 0,1\}$ (see Proposition \ref{P2.9}), we have
\begin{align*}
	(L^{p_{0}(\cdot)},\ell^{q_{0}})^{\alpha_{0}}_{\rm a} \cap (L^{p_{1}(\cdot)},\ell^{q_{1}})^{\alpha_{1}}_{\rm a} = \overline{\mathcal{C}_{\rm c}^{\infty}}.
\end{align*}
This shows that the space  $(L^{p_{0}(\cdot)},\ell^{q_{0}})^{\alpha_{0}}_{\rm a} \cap (L^{p_{1}(\cdot)},\ell^{q_{1}})^{\alpha_{1}}_{\rm a}$ is dense in both $(L^{p_{0}(\cdot)},\ell^{q_{0}})^{\alpha_{0}}_{\rm a}$ and $(L^{p_{1}(\cdot)},\ell^{q_{1}})^{\alpha_{1}}_{\rm a}$. Thus, from \eqref{eq4.14}  and Proposition \ref{prop 2.8}, we get
\begin{align*}
	\left[\mathcal{H}\left(p_{0}^{\prime}(\cdot),q_{0}^{\prime},\alpha_{0}^{\prime}\right), \mathcal{H}\left(p_{1}^{\prime}(\cdot),q_{1}^{\prime},\alpha_{1}^{\prime}\right)\right]^{\theta} = 
	\left(	\left[ (L^{p_{0}(\cdot)},\ell^{q_{0}})^{\alpha_{0}}_{\rm a}, (L^{p_{1}(\cdot)},\ell^{q_{1}})^{\alpha_{1}}_{\rm a} \right]_{\theta}\right)^{*}.
\end{align*}
Therefore, using Theorem \ref{Th4.4} and the equality \eqref{eq2.6}, successively, we obtain 
\begin{align} \label{eq4.15}
	\left[\mathcal{H}\left(p_{0}^{\prime}(\cdot),q_{0}^{\prime},\alpha_{0}^{\prime}\right), \mathcal{H}\left(p_{1}^{\prime}(\cdot),q_{1}^{\prime},\alpha_{1}^{\prime}\right)\right]^{\theta} = 
	\left[(L^{p(\cdot)},\ell^{q})^{\alpha}_{\rm a} \right]^{*} = \mathcal{H}\left(p^{\prime}(\cdot),q^{\prime},\alpha^{\prime}\right).
\end{align}
The combination of \eqref{eq4.13} and \eqref{eq4.15} ends the proof.
\end{proof}




\end{document}